 \documentclass[a4paper,11pt]{article}
 \usepackage[plainpages=false]{hyperref}
 \usepackage{amsfonts,amsmath,amssymb,amsthm}
 \usepackage{latexsym,lscape,rawfonts,mathrsfs}

 \usepackage[dvips]{color}
 \usepackage{multicol}

 \usepackage{natbib}

 \usepackage[all]{xy}
 \usepackage{eufrak}
 \usepackage{makeidx}
 \usepackage{graphicx,psfrag}

 \usepackage{array,tabularx}

 \usepackage{setspace}

 \newcommand{\D}[2]{\ensuremath{ \frac{\partial{#1}}{\partial{#2}}}}

 \newcommand{\R}{\ensuremath{\mathbb{R}}}
 \newcommand{\Z}{\ensuremath{\mathbb{Z}}}
 \newcommand{\CP}{\ensuremath{\mathbb{CP}}}
 \newcommand{\st}{\ensuremath{\sqrt{-1}}}

 \newcommand{\ddb}{\ensuremath{\partial \bar{\partial}}}

 \newcommand{\KRf}{K\"ahler Ricci flow\;}
 \newcommand{\KRF}{K\"ahler Ricci Flow\;}

 \DeclareMathOperator{\Vol}{Vol}
 \DeclareMathOperator{\diam}{diam}

 \newcommand{\Blow}[1]{\ensuremath{\mathbb{CP}^2 \# {#1}\overline{\mathbb{CP}}^2}}

 \newcommand{\norm}[2]{{ \ensuremath{\|} #1 \ensuremath{\|}}_{#2}}
 \newcommand{\snorm}[2]{{ \ensuremath{\left |} #1 \ensuremath{\right |}}_{#2}}

 \makeatletter
 \def\ExtendSymbol#1#2#3#4#5{\ext@arrow 0099{\arrowfill@#1#2#3}{#4}{#5}}
 
 \makeatother

 \definecolor{hao}{rgb}{1,0.5,0}
 \definecolor{miao}{cmyk}{0.5,0,0.2,0.2}
 \definecolor{qiao}{gray}{0.96}

 \newtheorem{corollary}{Corollary}[section]
 \newtheorem{proposition}{Proposition}[section]
 \newtheorem{lemma}{Lemma}[section]

 \newtheorem{theorem}{Theorem}[section]

 \newtheorem{remark}{Remark}[section]

\title{Remarks on \KRF}
\author{Xiuxiong Chen\footnote{Partially supported by a NSF grant.}\;,  Bing Wang}
\date{}
\begin{document}

 \maketitle

\begin{abstract}
   We show the convergence of \KRf directly if the
   $\alpha$-invariant of the canonical class is greater than $\frac{n}{n+1}$.
   Applying these convergence theorems,
   we can give a \KRf proof of Calabi conjecture on such Fano
   manifolds.
   In particular,  the existence of KE metrics on
   a lot of Fano surfaces can be proved by flow method.  Note that
   this geometric conclusion (based on the same assumption) was established
   earlier via elliptic method by G. Tian (cf. \cite{Tian87}, \cite{Tian90}
   and \cite{Tian97}).  However, a new proof based on \KRf should be still interesting in its own right.
  \end{abstract}

 \maketitle

 \section{Introduction}
   On a Fano manifold, the \KRf  was introduced as a possible means to search for
K\"ahler Einstein (KE) metrics.
   Following Yau's estimate (\cite{Yau78}),
    H. D. Cao (\cite{Cao85})
   proved that  the \KRf with smooth initial metric always exists globally.
   On a KE manifold, the first named author and Tian
   showed that \KRf converges exponentially fast toward the KE
   metric if the initial metric has positive bisectional
   curvature (cf. \cite{CT1}, \cite{CT2}).
  One important feature
  of these two papers is that the authors introduced a family of functionals
  to obtain a uniform  $C^0$-estimate of the evolved potential function
   $\varphi (t, \cdot).\; $
  Once the uniform  $C^{0}$-estimate is established, geometry of the
  evolved K\"ahler metrics is completely controlled and the flow will
  converge in the holomorphic category.\\

   Using his famous $\mu$-functional,  Perelman proved that
   scalar curvature, diameter and normalized Ricci
   potential are all uniformly bounded along \KRf (cf. \cite{PST}).
     In 2002, he announced that the \KRf will always converge
  to the KE metric on any KE manifold.  This result was  generalized to
  manifolds with K\"ahler Ricci solitons by Tian and Zhu (\cite{TZ}).   Based on these  fundamental estimates,
  if we assume that initial metric has positive bisectional
   curvature,  the first named author,  S. Sun and Tian (\cite{CTS})
   proved that the \KRf will converge to a KE metric automatically.
   Consequently, they give a  Ricci flow proof of Frankel conjecture
   which was initially proved by Siu-Yau (cf.\cite{SiY}) and Morri (cf.\cite{Mo}) independently.  Priori
  to \cite{CTS},
   partial progress was made in various works, e.g., \cite{Chen} and \cite{PSSW1}.  The study of \KRf is very intense after G. Perelman's fundamental estimates.     In this short note, we don't plan to analyze all these works  in depth. However,  we want to list a few
    references here for the convenience of readers:  ~\cite{Se1}, ~\cite{PSSW2},
   ~\cite{PSS}, ~\cite{Hei}, ~\cite{Ru1}, ~\cite{CH},~\cite{TZs},~\cite{RZZ},~\cite{FS}
  and references therein.\\

  On a K\"ahler manifold, Calabi conjectured that it admits a KE metric
  whenever its first Chern class $c_1$ has a definite sign or $c_1$ vanishes.
  He pointed out that the existence of KE metric is equivalent to the solvability of some
  Monge-Ampere equation which can be attacked by continuity method.
 This famous Calabi conjecture was proved by S. T. Yau in his
 celebrated work ~\cite{Yau78} when $c_1 < 0$ or $c_1 =0$.
 The case of $c_1 < 0$ was also obtained independently by T. Aubin
 (cf. \cite{Au}). For the case $c_1 > 0$ (Fano manifolds), situation is much
 more delicate. In \cite{Tian87}, \cite{TY} and \cite{Siu}, the
 existence of KE metric was proved on some special Fano manifolds. In the celebrated
 work \cite{Tian90}, Tian  finally proved the Calabi conjecture in dimension
 $2$.  He showed that a Fano surface $M$ admits a KE metric if and
 only if its automorphism group is reductive.
 In higher dimensional Fano manifolds,
  few general results are known. In \cite{Tian87}, Tian introduced the $\alpha$-invariant of a Fano manifold $M^n$
  and subsequently proved that $C^0$-estimate holds
 if  $\alpha$-invariant is greater than $\frac{n}{n+1}$.
 In \cite{Tian90} and \cite{Tian97},
 he showed that $C^0$-estimate can be obtained from the properness of
 $F$-functional.  Therefore, the Monge-Ampere equation is solvable
 whenever the $F$-functional is proper or the $\alpha$-invariant of
 $M^n$ is greater than $\frac{n}{n+1}$. Consequently, there is a KE metric in the
 canonical class of such Fano manifolds.
\\

    Inspired by these famous works in complex Monge-Ampere equation, we attempt to develop
    some estimates about the potential
    function $\varphi$ over K\"ahler Ricci flow.   The core issue is to obtain the $C^0$-estimate along the flow.
    Like the way used in continuity method, one can reduce the $C^0$-estimate of potential function $\varphi$
    to an integral estimate of $\varphi.\;$  This point was already
    observed by Yanir Rubinstein in~\cite{Ru1}.   As an easy
    application, we can prove directly  that potential function is uniformly bounded along the flow if
    $\alpha$-invariant is greater than $\frac{n}{n+1}$ or $F$ functional is proper.
    However, the $C^0$-estimate of $\varphi$ under the condition $F$-functional
    being proper  was proved  by Tian and Zhu~\cite{TZp} (cf. Proposition3.1 of \cite{TZ} also).
   We thank Tian for pointing this out to us.
   For the completeness of our presentation, we include a proof for the same statement. \\

   After we obtain the uniform $C^0$-norm of $\varphi$,  from the
   \KRf equation
   \begin{align}
     \dot{\varphi} = \log \frac{\omega_{\varphi}^n}{\omega^n} +
     \varphi - h_{\omega}.
     \label{eqn: krf}
   \end{align}
   and the free boundedness of $\dot{\varphi}$ (Lemma~\ref{lemma: dotphi}),
   one can easily see that every metric $\omega_{\varphi}$ is uniformly  equivalent to
   the initial metric $\omega$.  Moreover, the method in
   section 6 and 7 of \cite{CT2} applies directly on equation (\ref{eqn: krf}).  It follows that
   all $k$-th derivatives of $\varphi$ in the fixed gauge are uniformly bounded.
   Therefore, for every sequence $t_i \to \infty$, there is a subsequence
  $t_{i_k}$ such that
   $\omega+ \st \ddb \varphi(t_{i_k})$ smoothly converges to a limit metric $\omega + \ddb \varphi(\infty)$
   in the same gauge.
   This limit metric must be a KE metric since it is a K\"ahler Ricci Soliton metric
   with constant scalar curvature (See Section 4 for more details).    Therefore, the existence of
   KE metric is already proved.   Furthermore, by considering
   the first eigenvalue of $\triangle_{\omega_{\varphi}}$, one can even show that this
   flow converges to a unique KE metric exponentially
   fast(cf.~\cite{CT1}).\\

  As applications of our theorems, we can prove the existence of
  KE metrics on a lot of Fano surfaces by flow method.
  By classification theory,
  every Fano surface with zero Futaki invariant is either $\CP^1 \times
  \CP^1$, $\CP^2$ or diffeomorphic to $\Blow{k}, 3 \leq k \leq 8$.
  As $\CP^1 \times \CP^1$, $\CP^2$ and $\Blow{3}$ are all toric
  surfaces, the K\"ahler Ricci flows on them are studied in \cite{CW}.  The only
  interesting cases  are $4 \leq k \leq 8$.  Starting from any metric in canonical class,
  we can show the $C^0$-norms of potential functions are uniformly
  bounded if $M \sim \Blow{k}, k=4,5,7$.      For $M \sim \Blow{6}$,
  we can only show that $C^0$-estimate hold for some complex
  structure with nice symmetry.   If the symmetry of the complex
  structure is bad, we need to develop other methods to obtain the
  $C^0$-estimate, which will be discussed together with $k=8$ case
  in a subsequent paper.\\

  The organization of this paper is as follows. In section 2,
  we setup the notations and basic properties of K\"ahler geometry.
  In section 3, we write down some $C^0$-estimates of potential
  function $\varphi$ along K\"ahelr Ricci flow.  As applications of these estimates, we obtain some
  convergence theorems of \KRf and  apply them on Fano surfaces
  to prove Calabi conjecture in section 4.\\

  The authors would like to remark that the results of this paper
  have important overlaps with the results of~\cite{Ru1}.  Many theorems of
  this paper are actually implied in~\cite{Ru1}.  All of them are developed
  by the authors without  being aware of the results of [Ru1].
  We thank Y.A. Rubinstein for pointing these overlaps out to us.\\

 \noindent {\bf Acknowledgement:}  Both authors would like to thank G.Tian
 for many interesting and insightful discussions on the K\"ahler Ricci flow.
 The bulk of this work was carried out in department of Mathematics,
 University of Wisconsin at Madison.

 \section{Basic K\"ahler Geometry}
 Let $M$ be an $n$-dimensional compact K\"ahler manifold. A K\"ahler
 metric can be given by its K\"ahler form $\omega$ on $M$. In local
 coordinates $z_1, \cdots, z_n$, this $\omega$ is of the form
 \[
 \omega = \sqrt{-1} \displaystyle \sum_{i,j=1}^n\;g_{i \overline{j}}
 d\,z^i\wedge d\,z^{\overline{j}}  > 0,
 \]
 where $\{g_{i\overline {j}}\}$ is a positive definite Hermitian
 matrix function. The K\"ahler condition requires that $\omega$ is a
 closed positive (1,1)-form.  Given a K\"ahler metric $\omega$, its
 volume form  is
 \[
   \omega^n = {1\over {n!}}\;\left(\sqrt{-1} \right)^n \det\left(g_{i \overline{j}}\right)
  d\,z^1 \wedge d\,z^{\overline{1}}\wedge \cdots \wedge d\,z^n \wedge d\,z^{\overline{n}}.
 \]
 The curvature tensor is
 \[
  R_{i \overline{j} k \overline{l}} = - {{\partial^2 g_{i \overline{j}}} \over
 {\partial z^{k} \partial z^{\overline{l}}}} + \displaystyle
 \sum_{p,q=1}^n g^{p\overline{q}} {{\partial g_{i \overline{q}}}
 \over {\partial z^{k}}}  {{\partial g_{p \overline{j}}} \over
 {\partial z^{\overline{l}}}}, \qquad\forall\;i,j,k,l=1,2,\cdots n.
 \]
 The Ricci curvature form is
 \[
   {\rm Ric}(\omega) = \sqrt{-1} \displaystyle \sum_{i,j=1}^n \;R_{i \overline{j}}(\omega)
 d\,z^i\wedge d\,z^{\overline{j}} = -\sqrt{-1} \partial
 \overline{\partial} \log \;\det (g_{k \overline{l}}).
 \]
 It is a real, closed (1,1)-form and $[Ric]=2\pi c_1(M)$.\\

  From now on we we assume $M$ has positive first Chern class, i.e., $c_1(M)>0$.
  We call $[\omega]$ as a canonical K\"ahler class if $[\omega]=[Ric]=2\pi c_1(M)$.
    If we require the initial metric is in canonical class,
    then the normalized Ricci flow (c.f. \citet*{Cao85}) on $M$ is
\begin{equation}
  {{\partial g_{i \overline{j}}} \over {\partial t }} = g_{i \overline{j}}
  - R_{i \overline{j}}, \qquad\forall\; i,\; j= 1,2,\cdots,n.
\label{eq:kahlerricciflow}
\end{equation}
  Denote $\omega = \omega_{g(0)}$, $\omega_{g(t)} = \omega + \st \ddb \varphi_t$.
 $\varphi_t$ is called the K\"ahler potential and sometime it is denoted as
 $\varphi$ for simplicity.  On the level of K\"ahler potentials, \KRf becomes
\begin{equation}
   {{\partial \varphi} \over {\partial t }}
   =  \log {{\omega_{\varphi}}^n \over {\omega}^n } + \varphi -
   h_{\omega},
\label{eq:flowpotential}
\end{equation}
where $h_{\omega}$ is defined by
\begin{align*}
  {\rm Ric}(\omega)- \omega = \sqrt{-1} \partial \overline{\partial} h_{\omega},
    \; {\rm and}\;\displaystyle \int_M \;
  (e^{h_{\omega}} - 1)  {\omega}^n = 0.
\end{align*}
As usual, the flow equation (\ref{eq:kahlerricciflow}) or
(\ref{eq:flowpotential}) is referred as the \KRf in canonical class
of $M$. It is proved by Cao \citet{Cao85}, who followed Yau's
celebrated work \citet{Yau78}, that this flow exists globally for
any smooth initial K\"ahler metric in the canonical class.

   In this note, we only study \KRf in the canonical class. For the
   simplicity of notation, we may not mention that the flow is in
   canonical class every time.

 Let $\mathscr{P}(M, \omega)= \{\varphi | \omega + \st \ddb \varphi
 >0\}$. It is shown in~\cite{Tian87} that there is a small constant $\delta>0$ such that
 \begin{align*}
     \sup_{\varphi \in \mathscr{P}(M, \omega)}
    \frac{1}{V} \int_M e^{-\delta(\varphi -\sup_M \varphi)} \omega^n < \infty.
 \end{align*}
 The supreme of such $\delta$ is called the $\alpha$-invariant of $(M, \omega)$ and
 it is denoted as $\alpha(M, \omega)$.
 Let $G$ be a compact subgroup of $Aut(M)$ and $\omega$ is a $G$-invariant form. We denote
 \begin{align*}
   \mathscr{P}_G(M, \omega)
   =\{\varphi | \omega + \st \ddb \varphi >0, \varphi \; \textrm{is invariant under} \; G
   \}.
 \end{align*}
 Similarly, we can define $\alpha_G(M, \omega)$.

 For every $\varphi \in \mathscr{P}(M, \omega)$,
 we denote $\omega_{\varphi} = \omega + \st \ddb \varphi$.
 There are some well known functionals defined in this space.
 \begin{align}
    I_{\omega}(\varphi) & \triangleq \frac{1}{V} \int_M \varphi (\omega^n - \omega_{\varphi}^n) \notag\\
    &= \frac{1}{V}\sum_{i=0}^{n-1} \int_M \st \partial \varphi \wedge \bar{\partial}
          \varphi \wedge  \omega^i \wedge     \omega_{\varphi}^{n-1-i}. \label{eqn: I}\\
    J_{\omega}(\varphi) & \triangleq \frac{1}{V} \sum_{i=0}^{n-1}
    \frac{i+1}{n+1} \int_M \st \partial \varphi \wedge \bar{\partial}
    \varphi \wedge \omega^i \wedge \omega_{\varphi}^{n-1-i}. \label{eqn: J}\\
     F_{\omega}^0(\varphi) &\triangleq J_{\omega}(\varphi)
    - \frac{1}{V} \int_M \varphi \omega^n.  \label{eqn: F0}\\
    F_{\omega}(\varphi) &\triangleq  F_{\omega}^0(\varphi)- \log (\frac{1}{V} \int_M e^{h_{\omega} - \varphi}
    \omega^n). \label{eqn: F}\\
    \nu_{\omega(\varphi)} & \triangleq F_{\omega}(\varphi) +
    \frac{1}{V} \int_M h_{\omega} \omega^n - \frac{1}{V} \int_M
    h_{\omega_{\varphi}} \omega_{\varphi}^n.
    \label{eqn: nu}
 \end{align}
 The last functional is the well known Mabuchi K-energy.  It is generally defined  by its derivative.
 The formula here is proved in~\cite{DT}.

 We say $F_{\omega}$ is proper on $\mathscr{P}(M, \omega)$ if there
 exists an increasing function $\mu: \R \to [c, \infty)$ satisfying
 \begin{align*}
     \lim_{t \to \infty} \mu(t) = \infty; \quad
     F_{\omega}(\varphi) \geq \mu(I_{\omega}(\varphi)), \; \forall
     \varphi \in \mathscr{P}(M, \omega).
 \end{align*}
 Here $c$ is some number. In particular, $F_{\omega}$ is proper
 implies that $F_{\omega}$ is bounded from below.

 We list some basic properties of these functionals without giving
 proofs. Interested readers are referred to~\cite{Di}, ~\cite{Tb}, \cite{CT1} and references therein
 for more details.

 \begin{proposition}
  Suppose $M$ to be a Fano manifold without nontrivial holomorphic
  vector field.  If $M$ admits a KE metric in its canonical class
  $[\omega]$, then $F_{\omega}$
  is proper in $\mathscr{P}(M, \omega)$.
  \label{proposition: proper}
 \end{proposition}

 \begin{proposition}
     Along \KRf, $\D{}{t}F_{\omega}(\varphi_t) \leq 0$,
     $\D{}{t}\nu_{\omega}(\varphi_t) \leq 0$.
  \label{proposition: fmonotone}
 \end{proposition}

\section{Estimates along \KRf}

 As \KRf on Riemannian surface is very clear, we only consider
 \KRf on  Fano manifolds with complex dimension $n \geq 2$.

 We first list some well-known estimates.

  In his unpublished work, Perelman got some deep estimates along
 \KRf. We list his estimates below. The detailed proof can be found in Sesum and Tian's note ~\citet{PST}.
 \begin{proposition}[Perelman's Estimates]
 Suppose $\{(M^n, g(t)), 0 \leq t < \infty \}$ to be a \KRf solution. There
 are two positive constants $D, \kappa$ depending only on this flow such that the
 following two estimates hold.
 \begin{enumerate}
 \item  Let $R_{g_t}$ be the scalar curvature under metric $g_t$,
 $h_{\omega_{\varphi(t)}}$ be the Ricci potential of form $\omega_{\varphi(t)}$ satisfying $\frac{1}{V} \int_M e^{h_{\omega_{\varphi(t)}}} \omega_{\varphi(t)}^n=1$. Then we have
 \begin{align*}
      \norm{R_{g_t}}{C^0} + \diam_{g_t} M +
      \norm{h_{\omega_{\varphi(t)}}}{C^0} + \norm{\nabla h_{\omega_{\varphi(t)}}}{C^0} < D.
 \end{align*}
 \item   $ \displaystyle
      \frac{\Vol(B_{g_t}(x, r))}{r^{2n}} > \kappa$ for every
       $r \in (0, 1)$, $(x, t) \in M \times [0, \infty)$.
 \end{enumerate}
 \label{proposition: perelman}
 \end{proposition}

 In~\cite{Zhang} and~\cite{Ye}, Zhang and Ye obtained independently
 that Sobolev constant is uniformly bounded along every
 \KRf solution on a Fano manifold.

 \begin{proposition}[Sobolev Constant Estimate]
 There is a uniform Sobolev constant $C_S$ along the \KRf solution
 $\{(M^n, g(t)), 0 \leq t <\infty \}$. In other words, for every
 $f \in C^{\infty}(M)$, we have
 \begin{align*}
       (\int_M |f|^{\frac{2n}{n-1}}
       \omega_{\varphi}^n)^{\frac{n-1}{n}}
       \leq
       C_S(\int_M |\nabla f|^2 \omega_{\varphi}^n
       + \frac{1}{V^{\frac{1}{n}}} \int_M |f|^2  \omega_{\varphi}^n).
 \end{align*}
 \label{proposition: sobolev}
 \end{proposition}

 As an application of Perelman's estimate, weak Poincar\`e constant is also
 bounded along \KRf.

 \begin{proposition}[Weak Poincar\`e Constant Estimate]
   There is a uniform weak Poincar\`e constant $C_P$ along the \KRf
   solution $\{(M^n, g(t)), 0 \leq t < \infty \}$. In other words, for every
   nonnegative function $f \in C^{\infty}(M)$, we have
   \begin{align*}
        \frac{1}{V} \int_M f^2 \omega_{\varphi}^n   \leq
          C_P\{\frac{1}{V} \int_M |\nabla f|^2 \omega_{\varphi}^n
           + (\frac{1}{V} \int_M f  \omega_{\varphi}^n)^2\}.
   \end{align*}
 \label{proposition: poincare}
 \end{proposition}
 More details about this Proposition can be found in~\cite{TZ}.\\

 \begin{lemma}[c.f.~\cite{PSS}]
    By properly choosing initial condition, we have
 \begin{align*}
 \norm{\dot{\varphi}}{C^0} + \norm{\nabla \dot{\varphi}}{C^0}<C
 \end{align*}
 for some constant $C$ independent of time $t$.
 \label{lemma: dotphi}
 \end{lemma}

 From now on, we always choose initial condition properly such that
 $\dot{\varphi}$ is uniformly bounded.   Remember $\varphi$
 satisfies the equation
 $\dot{\varphi} = \log \frac{\omega_{\varphi}^n}{\omega^n} + \varphi - h_{\omega}$.
 This equation gives us a lot information.

  \begin{remark}
   When Sobolev constants, Poincar\`e constants and $\dot{\varphi}$ are all
   uniformly bounded, we can follow directly from continuity method
   to estimate $\snorm{\varphi}{C^0}$ from integrations of
   $\varphi$. This process is described from Lemma~\ref{lemma: supphi} to
   Theorem~\ref{theorem: conditions}.  However, from a conversation
   with Yanir Rubinstein, we know that this similarity was already observed by Yanir. Rubinstein~\cite{Ru1}.
   So these estimates were already implied in~\cite{Ru1} ( See~\cite{Ru2} also).  Just for the convenience of the readers,
   we include a complete proof here.
 \end{remark}

 \begin{lemma}
   There is a constant $C$ such that
 \begin{align}
   \frac{1}{V} \int_M (-\varphi) \omega_{\varphi}^n \leq
  n \sup_M \varphi  - \sum_{i=0}^{n-1}
     \frac{i}{V} \int_M \st \partial \varphi \wedge
     \bar{\partial} \varphi \wedge \omega^i \wedge
     \omega_{\varphi}^{n-1-i} + C.
  \label{lemma: dsupphi}
 \end{align}
 In particular, we have
   \begin{align}
       \frac{1}{V} \int_M (-\varphi) \omega_{\varphi}^n  \leq n \sup_M
       \varphi + C.
   \label{eqn: supphi}
   \end{align}
 \label{lemma: supphi}
 \end{lemma}

 \begin{proof}
   According to Proposition~\ref{proposition: fmonotone}, $F_{\omega}(\varphi)$ is
   non-increasing along \KRf.  Therefore, we have
   \begin{align*}
        0&=F_{\omega}(0) \\
        &\geq F_{\omega}(\varphi)\\
        &=J_{\omega}(\varphi) -\frac{1}{V} \int_M \varphi \omega^n
        -\log (\frac{1}{V} \int_M e^{h_{\omega} -\varphi}
        \omega^n)\\
        &=J_{\omega}(\varphi) - \frac{1}{V} \int_M \varphi \omega^n
        - \log(\frac{1}{V} \int_M e^{-\dot{\varphi}}
        \omega_{\varphi}^n)\\
        &=J_{\omega}(\varphi) -I_{\omega}(\varphi) + \frac{1}{V} \int_M (-\varphi) \omega_{\varphi}^n
                            - \log(\frac{1}{V} \int_M e^{-\dot{\varphi}}
                            \omega_{\varphi}^n).
   \end{align*}
  It follows that
  \begin{align*}
   \frac{1}{V} \int_M (-\varphi) \omega_{\varphi}^n
   &\leq I_{\omega}(\varphi) - J_{\omega}(\varphi)
    + \log(\frac{1}{V} \int_M e^{-\dot{\varphi}}
    \omega_{\varphi}^n).
  \end{align*}
  Plugging the expression of functional $I$ and $J$, we have
  \begin{align*}
   \frac{1}{V} \int_M (-\varphi) \omega_{\varphi}^n
   &\leq \sum_{i=0}^{n-1} \frac{n-i}{n+1} \frac{1}{V} \int_M
    \st \partial \varphi \wedge \bar{\partial} \varphi \wedge \omega^i
   \wedge \omega_{\varphi}^{n-1-i}+ \log(\frac{1}{V} \int_M e^{-\dot{\varphi}}   \omega_{\varphi}^n) \\
   &\leq  \frac{n}{n+1} I_{\omega}(\varphi)
     - \sum_{i=0}^{n-1} \frac{i}{n+1} \frac{1}{V} \int_M \st \partial \varphi \wedge
     \bar{\partial} \varphi \wedge \omega^i \wedge
     \omega_{\varphi}^{n-1-i}
     + \log(\frac{1}{V} \int_M e^{-\dot{\varphi}}   \omega_{\varphi}^n)\\
   &=\frac{n}{n+1} \frac{1}{V} \int_M \varphi \omega^n +
   \frac{n}{n+1} \frac{1}{V} \int_M (-\varphi) \omega_{\varphi}^n
   +\log(\frac{1}{V} \int_M e^{-\dot{\varphi}}
    \omega_{\varphi}^n)\\
   &\qquad - \sum_{i=0}^{n-1} \frac{i}{n+1} \frac{1}{V} \int_M \st \partial \varphi \wedge
     \bar{\partial} \varphi \wedge \omega^i \wedge
     \omega_{\varphi}^{n-1-i}.
  \end{align*}
  Note that $\dot{\varphi}$ is bounded, so we have
  \begin{align*}
   \frac{1}{V} \int_M (-\varphi) \omega_{\varphi}^n
   &\leq
   \frac{n}{V} \int_M \varphi \omega^n  - \sum_{i=0}^{n-1}
     \frac{i}{V} \int_M \st \partial \varphi \wedge
     \bar{\partial} \varphi \wedge \omega^i \wedge
     \omega_{\varphi}^{n-1-i} + C\\
   &\leq n \sup_M \varphi  - \sum_{i=0}^{n-1}
     \frac{i}{V} \int_M \st \partial \varphi \wedge
     \bar{\partial} \varphi \wedge \omega^i \wedge
     \omega_{\varphi}^{n-1-i} + C.
  \end{align*}
 \end{proof}

 \begin{lemma}
   There is
   a small constant $\delta$ and a big constant $C$ depending only on this flow such that
   \begin{align}
     \sup_M \varphi < \frac{1-\delta}{\delta}\int_M (-\varphi) \omega_{\varphi}^n  + C.
      \label{eqn: deltainv}
   \end{align}
 \label{lemma: nesup}
 \end{lemma}

 \begin{proof}

    As $M$ is a Fano manifold, we know there are two constants $\delta,
   C_{\delta}$ depending only on the initial metric $\omega$ such that
   \begin{align*}
         \sup_{\varphi \in \mathscr{P}(M, \omega)} \frac{1}{V}\int_M e^{-\delta(\varphi - \sup_M
         \varphi)} \omega^n < C_{\delta}.
   \end{align*}
   In particular, for every time $t$, $\varphi_t$ satisfies this inequality. Along \KRf, we have
    $\frac{\omega^n}{{\omega_{\varphi}}^n}= e^{\varphi -h_{\omega}-
    \dot{\varphi}}$. Plugging this into the previous inequality
    yields that
    \begin{align*}
       \frac{1}{V} \int_M e^{(1-\delta)\varphi + \delta \sup_M \varphi - h_{\omega}-
       \dot{\varphi}} \omega_{\varphi}^n < C_\delta.
    \end{align*}
    The convexity of exponential function tells us
    \begin{align*}
        \frac{1}{V} \int_M  \{\delta \sup_M \varphi +
        (1-\delta)\varphi - h_{\omega}- \dot{\varphi}\}
        \omega_{\varphi}^n < \log C_{\delta}.
    \end{align*}
   Therefore there is a constant $C$ depending only on this flow such that
   \begin{align}
     \sup_M \varphi < \frac{1-\delta}{\delta} \frac{1}{V} \int_M
     (-\varphi) \omega_{\varphi}^n  + C.
   \end{align}
  \end{proof}

 \begin{lemma}
  There is a constant $C$ depending only on the flow such that
  \begin{align*}
      \norm{\varphi}{C^0(M)} < C(\max\{0, \sup_M \varphi\} +1).
  \end{align*}
 \label{lemma: upperC0}
 \end{lemma}
 \begin{proof}
  Lemma~\ref{lemma: supphi} tells us
  $\displaystyle \frac{1}{V} \int_M (-\varphi) \omega_{\varphi}^n  \leq n \sup_M \varphi +
  C_1$.
  So we can choose constant $C_2 > n$ such that
  \footnote{$C_1, C_2, \cdots$ have different meaning in different Theorems.}
   \begin{align}
     \frac{1}{V} \int_M (-\varphi) \omega_{\varphi}^n  \leq  C_2 ( \max \{0, \sup_M \varphi \} + 1).
   \label{eqn: nphiin}
   \end{align}
   Clearly,
   $ \displaystyle
     \sup_M \varphi \leq \max \{0, \sup_M \varphi \}
      < C_2( \max \{0, \sup_M \varphi \} + 1).$
  Define
  \begin{align*}
  \bar{\varphi} \triangleq \varphi - 2C_2( \max \{0, \sup_M \varphi \} + 1).
  \end{align*}
  We have
  \begin{align*}
    \sup_M \bar{\varphi} < -1, \quad \frac{1}{V} \int_M
    |\bar{\varphi}| \omega_{\varphi}^n = \frac{1}{V} \int_M
    (-\bar{\varphi}) \omega_{\varphi}^n
    <3C_2( \max \{0, \sup_M \varphi \} +1).
  \end{align*}
   Direct computation shows that
  \begin{align}
      \int_M |\nabla |\bar{\varphi}|^{\frac{p+1}{2}}|^2
      \omega_{\varphi}^n &=
      \int_M \st \partial |\bar{\varphi}|^{\frac{p+1}{2}} \wedge \bar{\partial}
      |\bar{\varphi}|^{\frac{p+1}{2}} \wedge
      \omega_{\varphi}^{n-1} \notag\\
  &= \frac{(p+1)^2}{4} \int_M  \st |\bar{\varphi}|^{p-1} \partial
  |\bar{\varphi}| \wedge \bar{\partial} |\bar{\varphi}| \wedge
  \omega_{\varphi}^{n-1} \notag\\
  &=\frac{(p+1)^2}{4p} \int_M \st \partial |\bar{\varphi}|^p \wedge
  \bar{\partial} |\bar{\varphi}| \wedge
  \omega_{\varphi}^{n-1} \notag\\
  &=-\frac{(p+1)^2}{4p} \int_M  |\bar{\varphi}|^p  \st \partial \bar{\partial} |\bar{\varphi}| \wedge
  \omega_{\varphi}^{n-1} \notag\\
  &=\frac{(p+1)^2}{4p} \int_M  |\bar{\varphi}|^p  \st \partial \bar{\partial} \bar{\varphi} \wedge
  \omega_{\varphi}^{n-1} \notag\\
  &=\frac{(p+1)^2}{4p} \int_M  |\bar{\varphi}|^p  (\omega_{\varphi} - \omega) \wedge
  \omega_{\varphi}^{n-1} \notag\\
  & \leq \frac{(p+1)^2}{4p} \int_M  |\bar{\varphi}|^p
  \omega_{\varphi}^n.
  \label{eqn: phi12}
  \end{align}
  Let $p=1$, we have $\int_M |\nabla |\bar{\varphi}||^2 \omega_{\varphi}^n  \leq \int_M |\bar{\varphi}|
  \omega_{\varphi}^n$. Applying Poincar\`e  inequality (Proposition~\ref{proposition: poincare})  to
  $-\bar{\varphi}$ yields
   \begin{align*}
   \int_M |\bar{\varphi}|^2 \omega_{\varphi}^n &\leq
    C_P (\int_M |\nabla \bar{\varphi}|^2 \omega_{\varphi}^n
    + \frac{1}{V} (\int_M |\bar{\varphi}| \omega_{\varphi}^n)^2)\\
   &\leq C_P (\int_M |\bar{\varphi}| \omega_{\varphi}^n + \frac{1}{V} (\int_M |\bar{\varphi}|
   \omega_{\varphi}^n)^2).
  \end{align*}
  Using inequality (\ref{eqn: nphiin}), we obtain
  \begin{align}
    (\int_M |\bar{\varphi}|^2 \omega_{\varphi}^n)^{\frac12}
     \leq C_3 ( \max \{0, \sup_M \varphi \} + 1).
  \label{eqn: cfrompoincare}
  \end{align}
  For general $p \geq 1$, inequality (\ref{eqn: phi12}) can be rewritten as
  \begin{align*}
       \int_M |\nabla |\bar{\varphi}|^{\frac{p+1}{2}}|^2
      \omega_{\varphi}^n \leq \frac{(p+1)^2}{4p} \int_M  |\bar{\varphi}|^p
  \omega_{\varphi}^n < \frac{(p+1)^2}{4p} \int_M  |\bar{\varphi}|^{p+1} \omega_{\varphi}^n.
  \end{align*}
  Since Sobolev constants are uniformly bounded along \KRf,
  standard Moser iteration yields that
  \begin{align*}
    \norm{\bar{\varphi}}{C^0} &\leq C_4(\int_M |\bar{\varphi}|^2
    \omega_{\varphi}^n)^{\frac12}.
  \end{align*}
  Combining this with inequality (\ref{eqn: cfrompoincare}), we
  obtain
  \begin{align*}
    \norm{\bar{\varphi}}{C^0} \leq C_3 C_4 ( \max \{0, \sup_M \varphi \} + 1).
  \end{align*}
  Remember
   $\varphi = \bar{\varphi} + 2C_2( \max \{0, \sup_M \varphi \} + 1)$.
  Let $C= 2C_2 + C_3C_4 +1$, we have
  \begin{align*}
      \norm{\varphi}{C^0(M)} < C( \max \{0, \sup_M \varphi \} + 1).
  \end{align*}
 \end{proof}

 \begin{theorem}
  Along \KRf $\{(M^n, g(t)), 0 \leq t < \infty \}$ in the canonical class of
  Fano manifold $M$, the following conditions are equivalent.
 \begin{itemize}
 \item $\varphi$ is uniformly bounded.
 \item $\displaystyle \sup_M \varphi$ is uniformly bounded from above.
 \item $\displaystyle \inf_M \varphi$  is uniformly bounded from below.
 \item $\int_M \varphi \omega^n$ is uniformly bounded from above.
 \item $\int_M (-\varphi) \omega_{\varphi}^n$ is
   uniformly bounded from above.
 \item $I_{\omega}(\varphi)$ is uniformly bounded.
 \item $Osc_{M} \varphi$ is uniformly bounded.
 \end{itemize}
 \label{theorem: conditions}
 \end{theorem}

 \begin{proof}

 Look at Table~\ref{table: boundrelation}, it contains three
 circles: $(1234)$, $(256)$ and $(12789)$.   In order to prove this theorem, we only need
 to show the induction go through in every circle.   However, step 1
 is nothing but Lemma~\ref{lemma: nesup}, step 2 is just
 Lemma~\ref{lemma: upperC0}, steps 3, 4, 5, 7 are trivial.
 So only steps 6, 8, 9 need proof.

 \begin{table}
 \begin{displaymath}
   \xymatrix{& &\int_M (-\varphi) \omega_{\varphi}^n< C \ar@3{->}[d]^{1} & & I_{\omega}(\varphi)<C \ar@{~>}[ll]_{9}\\
    &\inf_M \varphi>-C \ar[ur]^{4} &\sup_M \varphi< C \ar@3{->}[d]^{2} & \int_M \varphi \omega^n <C \ar@{.>}[l]_{6} & \\
    & & \snorm{\varphi}{C^0(M)}< C \ar[ul]^{3} \ar@{.>}[ur]^{5} \ar@{~>}[rr]^{7}&
    & \ar@{~>}[uu]_{8} Osc_M \varphi < C}  \\
 \end{displaymath}
 \caption{The relations among bounds}
 \label{table: boundrelation}
 \end{table}

\noindent \textit{Step 6.   $\int_M \varphi \omega^n$
 bounded from above $\Rightarrow \displaystyle \sup_M \varphi$ bounded from above.}

  Since $\omega + \st \ddb \varphi >0$, take trace under metric $\omega$, we have $ -\triangle \varphi <
  n$. Plugging it into Green's function formula implies
  \begin{align*}
     \varphi(p) &= \frac{1}{V} \int_M \varphi \omega^n - \frac{1}{V}
     \int_M G(p, q) \triangle \varphi \omega^n(q)\\
      &\leq \frac{1}{V} \int_M \varphi \omega^n + \frac{n}{V} \int_M
      G(p, q)\omega^n(q)\\
      &\leq \frac{1}{V} \int_M \varphi \omega^n + C_1.
  \end{align*}
  Here $G$ is the nonnegative Green function of $(M, \omega)$, $p$
  is any point in $M$.  It follows that
  \begin{align*}
     \sup_M \varphi \leq \frac{1}{V} \int_M \varphi \omega^n + C_1.
  \end{align*}
  Therefore $\displaystyle \sup_M \varphi$ is bounded from above whenever
  $\int_M \varphi \omega^n$ is bounded from above.

\noindent \textit{Step 8. $Osc_M \varphi$ bounded $\Rightarrow
I_{\omega}(\varphi)$ bounded.}

 Note that $  I_{\omega}(\varphi)
    = \frac{1}{V} \int_M \varphi\omega^n - \frac{1}{V} \int_M  \varphi \omega_{\varphi}^n
    \leq \sup_M \varphi - \inf_M \varphi
    =Osc_M \varphi.$

\noindent \textit{Step 9. $I_{\omega}(\varphi)$ bounded
 $\Rightarrow \int_M (-\varphi) \omega_{\varphi}^n $ bounded from above.}

   The \KRf equation yields
   \begin{align}
      \frac{1}{V} \int_M (-\varphi) \omega^n &= \frac{1}{V} \int_M
      \{\log \frac{\omega_{\varphi}^n}{\omega^n} - \dot{\varphi} -h_{\omega} \}
      \omega^n \notag \\
      &\leq \frac{1}{V} \int_M \{-\dot{\varphi} - h_{\omega} \}
      \omega^n + \log \{ \frac{1}{V} \int_M \omega_{\varphi}^n \} \notag\\
      &= \frac{1}{V} \int_M \{-\dot{\varphi} - h_{\omega}\}   \omega^n \notag\\
      &< C_2.
      \label{eqn: phiorigin}
   \end{align}
   Here we used the fact that both $\dot{\varphi}$ and $h_{\omega}$ are
   uniformly bounded. Therefore, we have
   \begin{align}
        \frac{1}{V} \int_M (-\varphi) \omega_{\varphi}^n =
        I_{\omega}(\varphi) + \frac{1}{V} \int_M (-\varphi) \omega^n
         < I_{\omega}(\varphi) + C_2.
   \label{eqn: oscint}
   \end{align}
  This means that $I_{\omega}(\varphi)$ is bounded
    implies $\int_M (-\varphi) \omega_{\varphi}^n$ is bounded from above.
 \end{proof}

\section{Application of the Estimates}
   If $\varphi$ is uniformly bounded along the K\"ahler Ricci flow,
   then there must be
   a KE metric in the canonical class.  Actually, as
   we discussed in the introduction, there is a limit metric form
   $\omega_{\infty}=\omega + \st \ddb \varphi_{\infty}$ in the canonical class.
   As both $\varphi$ and $h_{\omega_{\varphi}}$ (Perelman's estimate)
   are uniformly bounded along the flow, it's easy to see from the
   definition of the functionals that $I_{\omega}(\varphi), J_{\omega}(\varphi), F_{\omega}(\varphi)$
   and $\nu_{\omega}(\varphi)$ are all bounded. In particular, the
   $K$-energy $\nu_{\omega}(\varphi)$ is bounded from  below. Therefore
   $\omega_{\infty}$ must be a metric with constant scalar
   curvature (cf. section 7 of \cite{CT2}, or ~\cite{PSSW1}). On the other hand, as a critical
   metric of Perelman's W-functional, $\omega_{\infty}$ must be a
   K\"ahler Ricci soliton.  So there is a smooth function $f$ such
   that
   \begin{align*}
      R_{i\bar{j}} + f_{i\bar{j}}- g_{i\bar{j}}=0, \quad
      f_{ij}=f_{i\bar{j}}=0.
   \end{align*}
   Taking trace yields $R + \triangle_{\omega_{\infty}} f -2=0$.
   Since $R$ is constant, $\triangle_{\omega_{\infty}} f$ has to be a constant
   and consequently zero. It follows that $f$ is a constant and we have
   $R_{i\bar{j}}=g_{i\bar{j}}$ under the metric $\omega_{\infty}$.
   This means $\omega_{\infty}$ is a KE metric in
   the canonical class.
   Notice  our convergence is in a fixed gauge, every limit KE metric form $\omega_{\infty}$ is
   compatible with the original complex structure.  There is no ``jump" of complex structure at all
   in this limit process.
   Using the method in \cite{CT1} and \cite{CT2}, we are able to show that the \KRf
   converges exponentially fast to a KE metric.

 Therefore, as corollary of Theorem~\ref{theorem: conditions}, we
 have the following theorems.

 \begin{theorem}
  $\{(M^n, g(t)), 0 \leq t < \infty \}$ is a \KRf solution
   initiating from a $G$-invariant metric $\omega$. $\alpha_G(M, \omega)> \frac{n}{n+1}$.
   Then this flow converges exponentially fast to a KE metric.
 \label{theorem: alpha}
 \end{theorem}

  \begin{proof}
   We only need to show $\varphi$ is uniformly bounded along the
   flow.

   Recall Lemma~\ref{lemma: supphi} and Lemma~\ref{lemma: nesup}.  Combining
   inequality (\ref{eqn: supphi}) and (\ref{eqn: deltainv}), we have
   \begin{align*}
      \frac{1}{V} \int_M (-\varphi) \omega_{\varphi}^n
      \leq n \sup_M \varphi
      + C_1 \leq n \cdot \frac{1-\delta}{\delta} \frac{1}{V} \int_M (-\varphi)
      \omega_{\varphi} + C_2.
   \end{align*}
   Since $\alpha_G> \frac{n}{n+1}$, we can choose $\delta >
   \frac{n}{n+1}$ such that  $n \cdot \frac{1-\delta}{\delta} < 1$. Therefore,
 \begin{align*}
      (1- n \frac{1-\delta}{\delta}) \frac{1}{V} \int_M (-\varphi)
      \omega_{\varphi}^n < C_2.
 \end{align*}
  It follows that both $\frac{1}{V} \int_M (-\varphi)
  \omega_{\varphi}^n$ and $\displaystyle \sup_M \varphi$ are uniformly bounded.
  So Theorem~\ref{theorem: conditions} implies that $\varphi$ is
  uniformly bounded.
 \end{proof}

 \begin{remark}
    This theorem implies the existence of KE metric on a lot of
  Fano manifolds. For example, every Fano manifold $M$ without
  $G$-invariant multiplier ideal sheaf  (See~\cite{Na})
  has $\alpha_G(M, \omega) \geq
  1>\frac{n}{n+1}$.  Therefore KE metric exists in its canonical
  class.

   On a Mukai-Umemura 3-fold $M$, Donaldson (\cite{Don}) showed that
    $\alpha_{SO(3)}(M, \omega)=\frac56>\frac34$.  Therefore
   Calabi conjecture holds on this manifold.

 \end{remark}

  \begin{theorem}[\cite{TZp}]
   $\{(M^n, g(t)), 0 \leq t < \infty \}$ is a \KRf solution initiating from $\omega$. $F_{\omega}$ is proper
    on the space  $\mathscr{P}(M^n, \omega)$.   Then this flow  converges exponentially fast to a KE metric.
 \label{theorem: proper}
 \end{theorem}

 \begin{proof}
    As $F$ is proper and
  $F_{\omega}(\varphi) \leq F_{\omega}(0)$ along \KRf, we see that  $I_{\omega}(\varphi)$
  is uniformly bounded. So $\varphi$ is uniformly bounded.
 \end{proof}

  Perelman has claimed that the \KRf will converge to the
  KE metric if KE metric exists in the
  canonical class (A generalization of this claim is proved in~\cite{TZ}).
  These two theorems can be achieved directly from
  Tian's existence theorem of KE metrics and
  Perelman's claim.    However, we obtain these theorems from
  Ricci flow without assuming the existence of KE metric.

 \begin{lemma}
    Suppose $(M^n, \omega)$ is a Fano manifold, $(N^n, \omega_{KE})$ is a K\"ahler Einstein
    manifold, $\pi: M \to N$ is a branched covering map satisfying  $[\pi^* \omega_{KE} ]=\lambda c_1(M)$
    for some number $\lambda>1$,  $G$ is the deck transformation group and $\omega$ is
    $G$-invariant.

    If  $\int_M (\frac{\pi^* \omega_{KE}^n}{\omega^n})^{-\gamma} \omega^n< \infty$ for some  $\gamma> \frac{1}{\lambda-1}$,
    then the \KRf initiating from $(M,  \omega)$ will converge
    exponentially fast to a KE metric. In particular,
    there exists a KE metric on $M$.
   \label{lemma: coverproper}
  \end{lemma}

 \begin{proof}
    For simplicity, we suppose this covering is $p$-sheeted, i.e.,
   $|G|=p$. We denote $f$ as the smooth function  $\frac{\pi^* \omega_{KE}^n}{\omega^n}$.
   Let $\varphi=\varphi_t$,
    $\pi^*(\frac{1}{\lambda} \omega_{KE})= \omega+ \st \ddb u$,
   we have
 \begin{align*}
     F_{\omega}(\varphi)&= F_{\omega}^0 (\varphi) - \log(\frac{1}{V} \int_M e^{h_{\omega}-\varphi}
     \omega^n) \notag \\
   &= F_{\omega}^0(u) + F_{\pi^*(\frac{1}{\lambda}
   \omega_{KE})}^0(\varphi-u) - \log(\frac{1}{V} \int_M e^{-\dot{\varphi}} \omega_{\varphi}^n). \notag\\
 \end{align*}
 Note that $u$ is a fixed function and $F_{\omega}^0(u)$ is a fixed
 number. So we have
 \begin{align}
   F_{\omega}(\varphi) & \geq F_{\pi^*(\frac{1}{\lambda} \omega_{KE})}^0(\varphi-u) - C \notag \\
   &=\frac{1}{\lambda} F_{\pi^*(\omega_{KE})}^0(\lambda(\varphi-u))-C \notag\\
   &=\frac{1}{\lambda} F_{\omega_{KE}}^0(\lambda \pi_*(\varphi-u))-C.
   \label{eqn: updown}
 \end{align}
 The last step is well defined since $\varphi$ is $G$-invariant.
 Notice that  $F_{\omega_{KE}}(\lambda \pi_*(\varphi-u)) $  is bounded from below.
 It follows from inequality (\ref{eqn: updown}) that
 \begin{align}
    F_{\omega}(\varphi) &\geq \frac{1}{\lambda}\{ F_{\omega_{KE}}(\lambda \pi_*(\varphi-u))
      + \log (\frac{1}{\Vol(N)} \int_N e^{\lambda \pi_*(u-\varphi)} \omega_{KE}^n)\}
    -C \notag \\
    &\geq \frac{1}{\lambda} \log (\frac{1}{p} \frac{1}{\Vol(N)} \int_M e^{\lambda(u-\varphi)}
    \pi^* \omega_{KE}^n) -C \notag\\
    &\geq \frac{1}{\lambda} \log ( \int_M e^{-\lambda \varphi} \pi^* \omega_{KE}^n) -C.
 \label{eqn: FE}
 \end{align}
 In the last step, we used the property that $u$ is bounded on $M$.

 Let $\beta= \frac{\lambda\gamma}{\gamma+1}>1$. H\"older inequality implies
 \begin{align*}
    \int_M e^{-\beta \varphi} \omega^n
   &=\int_M e^{-\beta \varphi} \cdot f^{\frac{\beta}{\lambda}} \cdot f^{-\frac{\beta}{\lambda}}
    \omega^n\\
 &\leq ( \int_M e^{-\lambda\varphi} f \omega^n)^{\frac{\gamma}{\gamma+1}}
 (\int_M f^{-\gamma} \omega^n)^{\frac{1}{\gamma+1}}.
 \end{align*}
 As $\int_M f^{-\gamma} \omega^n$ is a finite number, we
 have
 \begin{align*}
  \int_M e^{-\beta \varphi} \omega^n
  \leq  C (\int_M e^{-\lambda\varphi} \pi^* \omega_{KE}^n)^{\frac{\gamma}{\gamma+1}}.
 \end{align*}
 It follows that
 \begin{align*}
    \frac{\gamma}{\gamma+1} \log( \int_M e^{-\lambda\varphi} \pi^* \omega_{KE}^n)
    \geq - C + \log(\int_M e^{-\beta \varphi} \omega^n)
 \end{align*}
 Putting this inequality into (\ref{eqn: FE}) gives us
 \begin{align*}
  F_{\omega}(\varphi) &\geq \frac{(\gamma+1)}{\lambda\gamma}\log(\int_M e^{-\beta \varphi}
  \omega^n) -C\\
   &=\frac{(\gamma+1)}{\lambda\gamma} \log(\int_M e^{-\beta \varphi}
   \cdot e^{\varphi - h_{\omega} - \dot{\varphi}} \omega_{\varphi}^n) -C\\
   &\geq \frac{1}{\beta} \log(\int_M e^{(1-\beta)\varphi}
   \omega_{\varphi}^n) -C.
 \end{align*}
 The convexity of exponential map together with monotonicity of $F_{\omega}(\varphi)$ implies that
 \begin{align*}
 0= F_{\omega}(0) \geq  F_{\omega}(\varphi) \geq \frac{(\beta-1)}{\beta} \int_M
     (-\varphi) \omega_{\varphi}^n - C.
 \end{align*}
 As $\beta-1>0$, the previous inequality implies $\int_M (-\varphi)\omega_{\varphi}^n$
 is uniformly bounded from above along the flow. By
 Theorem~\ref{theorem: conditions}, we know $\varphi$ is uniformly
 bounded along the flow. Therefore this flow converges exponentially
 fast to a KE metric.

 \end{proof}

   If we denote $R(\pi) \subset M$ as the ramification divisor of $\pi$.
   Choose $x \in R(\pi)$,  let $s$ be the defining holomorphic
   function (locally defined ) of $R(\pi)$ at $x$.  Define
   \begin{align*}
      \alpha_x(R(\pi))& \triangleq \sup\{ \lambda \geq 0: |s|^{-2\lambda} \; \textrm{is } \; L^1 \; \textrm{on a neighborhood of
      x}.\}\\
        \alpha(R(\pi)) & \triangleq \inf_{x \in R(\pi)}
        \alpha_x(R(\pi)).
   \end{align*}
   Note that $\alpha(R(\pi))=1$ if $R(\pi)$ is a reduced smooth divisor.
   Denote
   \begin{align*}
       c \triangleq \sup \{\lambda \geq 0 |  \int_M (\frac{\pi^*\omega_{KE}^n}{\omega^n})^{-\lambda} \omega^n < \infty
       \}.
   \end{align*}
   It's shown in \cite{Ar} (Lemma 2.8) that $c= \alpha(R(\pi))$.  Therefore, we
   have the following theorem.

  \begin{theorem}
    Suppose $(M^n, \omega)$ is a Fano manifold, $(N^n, \omega_{KE})$ is a K\"ahler Einstein
    manifold, $\pi: M \to N$ is a branched covering map satisfying  $[\pi^* \omega_{KE} ]=\lambda c_1(M)$
    for some number $\lambda>1$. $R(\pi)$ is the ramification divisor of $\pi$,
    $G$ is the deck transformation group and $\omega$ is
    $G$-invariant.

    If  $\alpha(R(\pi))> \frac{1}{\lambda-1}$,
    then the \KRf initiating from $(M,  \omega)$ will converge
    exponentially fast to a KE metric. In particular,
    there exists a KE metric on $M$.
   \label{theorem: coverproper}
  \end{theorem}

  As in~\cite{Ar}, we can generalize this theorem to multiple
  covers. We first fix some notations.

  Let $D_0, \cdots D_l$ be divisors of $M$.
  Fix $\displaystyle x \in \bigcup_{i=0}^l D_i$. Let $f_i$ be the local defining holomorphic functions of
  $D_i$.   Define
  \begin{align*}
    \alpha_x(D_0, \cdots, D_l) &\triangleq \sup \{ \delta \geq 0:   (|f_0|+ \cdots |f_l|)^{-2
    \delta} \;
     \textrm{is} \; L^1 \; \textrm{on a neighborhood of} \; x \}.\\
    \alpha(D_0, \cdots, D_l) &\triangleq \inf_{x \in \bigcup_{i=0}^l
    D_i} \alpha_x(D_0, \cdots, D_l).
  \end{align*}
  Then Theorem~\ref{theorem: coverproper} can be generalized as the
  following.

  \begin{theorem}
    Suppose $(M^n, \omega)$ is a Fano manifold, $(N_i^n, \omega_{KE}),  0 \leq i \leq l$ is a K\"ahler Einstein
    manifold, $\pi_i: M \to N_i$ is a branched covering map satisfying  $[\pi^* \omega_{KE} ]=\lambda_i c_1(M)$
    for some number $\lambda_i>1$. $R(\pi_i)$ is the ramification divisor of $\pi_i$, $G_i$ is the deck transformation group
     and $\omega$ is
    $G_i$-invariant for every $0 \leq i \leq l$.

    If $\displaystyle \alpha(R(\pi_0), \cdots, R(\pi_l))> \max_{0 \leq i \leq l} \frac{1}{\lambda_i-1}$,
    then the \KRf initiating from $(M,  \omega)$ will converge
    exponentially fast to a KE metric.
    In particular, if $\displaystyle \bigcap_{i=0}^l R(\pi_i) =
     \emptyset$,  then the \KRf initiating from $(M,  \omega)$ will converge
    exponentially fast to a KE metric.
   \label{theorem: coverproperg}
  \end{theorem}

 Then let's apply previous theorems on Fano surfaces.

 \begin{corollary}[\cite{Hei}]
    Let $M$ be a Fano surface and  $M \sim \Blow{4}$.  $\omega$ is a
    metric form  in canonical class.
    Then the \KRf initiating  from $\omega$ converges exponentially
    fast to a KE metric.
 \end{corollary}

 \begin{proof}
   First we show that there exists a KE metric in the
   canonical class.  Start \KRf from a metric $\omega_0$ which is invariant under
   the finite automorphism group $G$. Remember $\alpha_G(M, \omega_0) > \frac23$ (See~\cite{TY}),
   Theorem~\ref{theorem: alpha} implies the existence of
   KE metric in the canonical class.

   Now we consider the \KRf initiating from any metric $\omega$ in the
   canonical class.   Note that $M$ has no nontrivial holomorphic vector
   field. By Proposition~\ref{proposition: proper},
   $F_{\omega}$ is proper since the existence of a KE metric.
   By Theorem~\ref{theorem: proper}, we obtain the result we need.
 \end{proof}

  \begin{corollary}
    Let $M$ be a Fano surface and  $M \sim \Blow{k}, \; k=5,7$.  $\omega$ is a
    metric form  in canonical class.
    Then the \KRf initiating  from $\omega$ converges exponentially
    fast to a KE metric.
 \end{corollary}
 \begin{proof}
   As $Aut(M)$ is discrete for each such $M$, the existence of KE
   metric will imply the properness of $F_{\omega}$ by
   Proposition~\ref{proposition: proper}. Therefore, by theorem~\ref{theorem: proper}, we only need to show the
   existence of a KE metric.   We will use \KRf with symmetry to find the
   KE metric.  The main tools are Theorem~\ref{theorem: coverproper}
   and Theorem~\ref{theorem: coverproperg}.  For the construction of
   the branched coverings, see~\cite{De} for details.

 \noindent
 \textit{Case1.  $M \sim \Blow{5}$.}

   Let $N_0=N_1=N_2=N_3=N_4=\CP^1 \times \CP^1$, we can
   construct $\pi_i: M \to N_i, 0 \leq i \leq 4$ such that
   \begin{align*}
     \bigcap_{i=0}^4 R(\pi_i) =\emptyset.
   \end{align*}
   Moreover, let $G_i$ be the deck transformation for covering
   $\pi_i$, we can find a metric $\omega$ which is invariant under
   every $G_i$.  Theorem~\ref{theorem: coverproperg} applies and  \KRf
   tends to a KE  metric in $\displaystyle \bigcap_{i=0}^4 \mathscr{P}_{G_i}(M, \omega)$.

   Actually,  $M \sim \Blow{5}$ can be embedded into $\CP^4$ as the complete intersection of
    two quadrics $Q_1$ and $Q_2$. By the result of Miles Reid (\cite{Re}), we can
    find a coordinate system of $\CP^4$ such that
    \begin{align*}
      Q_1&=\{ x_0^2+x_1^2+x_2^2+x_3^2 + x_4^2=0\}\\
      Q_2&=\{ a_0 x_0^2 + \cdots +a_4 x_4^2=0\}
    \end{align*}
  where $a_i \neq a_j$ when $i \neq j$.    As $M= Q_1 \cap Q_2$, by ``forgetting" $x_0$, we obtain a projection
  map $\pi_0: M \to  N_0 \subset \CP^3$ where
  \begin{align*}
     N_0= \{(x_1: x_2: x_3: x_4)| (a_1-a_0)x_1^2 + \cdots   (a_4-a_0)x_4^2=0\}.
  \end{align*}
  So $N_0$ is a smooth quadratic surface in $\CP^3$. It is
  biholomorphic to $\CP^1 \times \CP^1$ and it admits a KE metric.
  Similarly, we can define $\pi_i, 1 \leq i \leq 4$.   We have
  \begin{align*}
     \bigcap_{i=0}^4 R(\pi_i) \subset \bigcap_{i=0}^4 \{x_i=0\}= \emptyset.
  \end{align*}

   $G_i=\Z_2$ and it acts on $M$ by multiplying $\pm 1$ on the $i$-th coordinate.
   Clearly, by taking average, we can find a metric $\omega$ which is invariant under
   $\displaystyle \oplus_{i=0}^4 G_i$.

 \noindent
 \textit{Case2.  $M \sim \Blow{7}$.}

   In this case,  $M$ is a branched double cover of $N=\CP^2$.
   $R(\pi)$ is a smooth curve, so $\alpha(R(\pi))=1$.   $[\pi_*(\omega_{KE})]=3c_1(M)$, so $\lambda=3$.
   Clearly, $\alpha(R(\pi))=1 > \frac{1}{\lambda-1}=\frac12$.
   So Theorem~\ref{theorem: coverproper} applies and there exists a KE
   metric in  $\mathscr{P}_G(M, \omega)$.

 \end{proof}

 \begin{corollary}
    Let $M$ be a Fano surface and  $M \sim \Blow{6}$. It is well
    known that $M$ is a cubic surface in $\CP^3$.   Suppose that
    \begin{align*}
       M= \{x_0^3+ \cdots x_{l}^3 + f(x_l, \cdots, x_3)=0\}
       \subset \CP^3
    \end{align*}
    for some $l \geq 1$.    $\omega$ is a  metric form  in canonical class.
    Then the \KRf initiating  from $\omega$ converges exponentially
    fast to a KE metric.
 \label{corollary: cubic}
 \end{corollary}

 \begin{proof}
   As in the previous two corollaries, we only need to prove the
   existence of KE metric.

   In this case, $M$ admits $l+1$ branched covering $\pi_i: M \to
   N_i=\CP^2$ obtained by
 \begin{align*}
   \pi(x_0, \cdots, x_3) = (x_0, \cdots, \hat{x_i}, \cdots, x_3), \quad 0 \leq i \leq l.
 \end{align*}
  $G_i= \Z_3$ acts by multiplication of roots of $z^3=1$ on the $i$-th coordinate of $\CP^3$.
  Direct computation shows $[\pi_i^*(\omega_{FS})]= 3c_1(M)$ for
  every $0 \leq i \leq l$.  As argued in Proposition 3.1 of \cite{Ar}, we can obtain
 \begin{align*}
  \alpha(R(\pi_0), \cdots, R(\pi_{l}))> \frac{1}{3-1}=\frac12.
 \end{align*}
 So Theorem~\ref{theorem: coverproperg} applies.  Starting from an
 $\omega$ which is invariant under $\displaystyle
 \oplus_{i=0}^{l}G_i$, \KRf will converge to a KE metric.
 \end{proof}

  For those cubic surfaces ($\Blow{6}$) with bad symmetry, i.e., those cubic surfaces whose equations cannot be
  written as in Corollary~\ref{corollary: cubic}, we have
  other methods to study the behavior of \KRf on it, which will be
  discussed together with \KRf on $\Blow{8}$ in a subsequent
  paper.\\

\begin{remark}
   On $\Blow{4}$ and $\Blow{5}$ with
   special complex structure,  Gordon Heier (\cite{Hei}) proved the convergence
   of \KRf by multiplier ideal sheaf method.   This method is first
   studied in \cite{PSS} for general \KRf on Fano manifolds. It is
   improved by Yanir A. Rubinstein in \cite{Ru1}.
\end{remark}

 \vspace{0.5in}

 Xiuxiong Chen,  Department of Mathematics, University of
 Wisconsin-Madison, Madison, WI 53706, USA; xiu@math.wisc.edu\\

 Bing  Wang, Department of Mathematics, University of Wisconsin-Madison,
 Madison, WI, 53706, USA; bwang@math.wisc.edu

 \qquad \qquad \quad \; Department of Mathematics, Princeton University,
  Princeton, NJ 08544, USA; bingw@math.princeton.edu
\end{document}